\documentclass[a4paper,11pt]{article} 
\usepackage{amsmath,amssymb,amsfonts} 
\usepackage[all]{xy}
\newcommand{\comment}[1]{}

\DeclareMathOperator{\rank}{rank}
\DeclareMathOperator{\Pic}{Pic}

\renewcommand{\AA}{\mathbb{A}}
\newcommand{\CC}{\mathbb{C}}

\newcommand{\PP}{\mathbb{P}}

\newcommand{\ZZ}{\mathbb{Z}}

\newcommand{\al}{\alpha}
\newcommand{\be}{\beta}

\newcommand{\del}{\delta}
\newcommand{\eps}{\varepsilon}
\newcommand{\fie}{\varphi}
\newcommand{\ga}{\gamma}
\newcommand{\la}{\lambda}
\newcommand{\si}{\sigma}
\newcommand{\Fie}{\Phi}
\newcommand{\Gam}{\Gamma}
\newcommand{\Del}{\Delta}
\newcommand{\La}{\Lambda}
\newcommand{\Si}{\Sigma}

\newcommand{\coset}[1]{\overline{#1}}
\newcommand{\xprime}{x^\prime}

\newcommand{\eL}{\mathcal{L}}
\newcommand{\eM}{\mathcal{M}}

\newcommand{\Oh}{\mathcal{O}}

\newcommand{\Hom}{\text{Hom}}
\newcommand{\Ima}{\text{Im}}
\newcommand{\Ker}{\text{Ker}}
\newcommand{\End}{\text{End}}

\DeclareMathOperator{\Alb}{Alb}
\DeclareMathOperator{\Tors}{Tors}
\DeclareMathOperator{\alb}{alb}
\DeclareMathOperator{\inv}{inv}

\newcommand{\Psihat}{\widetilde\Psi}
\newcommand{\psihat}{\widetilde\psi}
\newcommand{\Xhat}{\widehat X}
\newcommand{\Xtilde}{\widetilde X}
\newcommand{\Shat}{\widehat S}

\newcommand{\isom}{\cong}

\newcommand{\Span}[1]{\left<#1\right>}
\newcommand{\paren}[1]{\left(#1\right)}

\newtheorem{thm}{Theorem}[section]

\newtheorem{lemma}[thm]{Lemma}
\newtheorem{prop}[thm]{Proposition}

\newtheorem{rmk}[thm]{Remark}
\newtheorem{dfn}[thm]{Definition}
\newenvironment{pf}{\paragraph{Proof}}{\par\medskip}
\newenvironment{pfthm}{\paragraph{Proof of theorem}}{\par\medskip}

\newenvironment{pflemma}{\paragraph{Proof of lemma}}{\par\medskip}
\newcommand{\qed}{\hfill\square}

\title{Kulikov surfaces form a connected component\\ of the moduli space}
\author{Tsz On Mario Chan, Stephen Coughlan}
\date{}

\begin{document}
\maketitle
\begin{abstract}
We show that the Kulikov surfaces form a connected component of the moduli space of surfaces of general type with $p_g=0$ and $K^2=6$. We also give a new description for these surfaces, extending ideas of Inoue. Finally we calculate the bicanonical degree of Kulikov surfaces, and prove that they verify the Bloch conjecture.
\end{abstract}
\section{Introduction}
In 1966, Burniat \cite{Bur} constructed some new examples of surfaces of general type with $p_g=0$ as $(\ZZ/2)^2$-covers of the plane branched in a certain configuration of lines. Then in 1994, Inoue \cite{Inoue} constructed examples with the same numerical invariants as Burniat's examples by taking a $(\ZZ/2)^3$-quotient of a hypersurface inside a product of three elliptic curves. It seemed to be well known that the two constructions were equivalent, and a proof was written down in \cite{BC}. The Inoue construction was then used in \cite{BC} to explicitly calculate the fundamental group and homology of the Burniat--Inoue surfaces. This enabled a proof that primary Burniat--Inoue surfaces (those with $K^2=6$) form a connected component of the moduli space of surfaces of general type, and further, this component is actually closed under homotopy equivalence (see below for an explanation of this terminology). The first proof that Burniat--Inoue surfaces form a connected component of the moduli space appeared in \cite{MP2}.

Kulikov \cite{Kul} exhibited a new example of a surface of general type with $p_g=0$ and $K^2=6$ as a $(\ZZ/3)^2$-cover of the plane branched in a different configuration of lines (see Figure \ref{fig!configs}). In this paper we show how the Kulikov surfaces fit into the same framework as Burniat--Inoue surfaces by exhibiting an Inoue-type construction. We go on to prove that Kulikov surfaces form a $1$-dimensional connected component of the moduli space, which is closed under homotopy equivalence.

The layout and main results of this paper are as follows: in Section \ref{sec!surfaces} we discuss preliminary results on abelian covers of the plane branched in line configurations, and use these to define the Kulikov surface. In Section \ref{sec!moduli} we prove that there is an Inoue-type construction for Kulikov surfaces, and that they form a $1$-dimensional irreducible component of the moduli space of surfaces of general type. We also describe three interesting degenerate Kulikov surfaces, which can be viewed as the boundary points of the compactified moduli space.

In Section \ref{sec!homotopy}, following ideas employed in \cite{BC}, we use the Inoue-type construction to explicitly calculate the fundamental groups and homology of our surfaces:
\begin{thm}
If $X$ is a Kulikov surface, then $\pi_1(X)=\Gam$ and $H_1(X,\ZZ)=(\ZZ/3)^3$.
\end{thm}
Here $\Gam$ is a certain infinite subgroup of the affine group $\AA(3,\CC)$, decribed in terms of generators in Section \ref{sec!homotopy}. Having determined the fundamental group, we are able to prove our second main result:
\begin{thm} The moduli space of Kulikov surfaces is closed under homotopy equivalence.
\end{thm}
In other words, any compact complex surface which is homotopy equivalent to a Kulikov surface is given by our construction. We say that the moduli space is  \emph{closed under homotopy equivalence}. This implies the pencil of Kulikov surfaces forms a connected component of the moduli space of surfaces of general type with $p_g=0$ and $K^2=6$. One might conjecture that any surface with the same $\pi_1$ is given by our construction. We have not been able to prove this, because the proof of Proposition \ref{prop!albanese}$(2)$ relies on homotopy equivalence.

Section \ref{sec!bicanonical} contains a calculation of the degree of the bicanonical map of the Kulikov surfaces:
\begin{thm} The bicanonical morphism of a Kulikov surface is birational.
\end{thm}
This is the expected result; by \cite{MP}, \cite{MP2} the possible values are $1$ or $2$, but it would be surprising if the degree did not divide $3$.
In the course of calculating the degree of the bicanonical map, we prove Prop\-osition \ref{prop!bican}, which gives a formula for the eigenspace decomposition of the bicanonical sheaf $\omega_X^2$ for an abelian cover. This complements the standard formulas of \cite{Par} Proposition 4.1, and we believe it is of independent interest.

In our last section, we verify the Bloch conjecture on zero cycles for Kulikov surfaces. Let $A_0^0(X)$ denote the group of zero cycles of degree $0$ on $X$, modulo rational equivalence.
\begin{thm} The Kulikov surface $X$ has $A_0^0(X)=0$.
\end{thm}
The theorem is proved using a well known method of \cite{IM}, which involves a careful examination of various quotient surfaces associated to the maximal abelian cover $\Xhat\to X$. Unfortunately the final step in our proof is done by computer algebra.

During the final preparation of this manuscript we were informed by V.~Alexeev that he and R.~Pardini also computed the degenerate Kulikov surfaces of Section \ref{sec!degenerate}, using similar methods to those of \cite{AP}.

\subsection*{Acknowledgements}
Both authors thank Ingrid Bauer and Fabrizio Catanese for their help and encouragement with this paper. We also thank both referees for their comments which improved our manuscript in various ways. The first author was supported by the DAAD Forschungsstipendien f\"ur Doktoranden and the second author was supported by DFG Forschergruppe $790$ ``Classification of algebraic surfaces and compact complex manifolds''.
\section{Surfaces as abelian covers}\label{sec!surfaces}
\subsection{Abelian covers}
We briefly review the preliminary results needed for our calculations. More details can be found in several articles, including \cite{BC2}, \cite{Cat}, \cite{Par}. Assume $Y$ is nonsingular and let $\fie\colon X\to Y$ be a finite morphism with a faithful action of the abelian group $G$ on $X$ so that $Y$ is the quotient. Let $\Delta$ be the branch divisor of $\fie$, then $\fie$ is determined by the surjective homomorphism $\Fie\colon\pi_1(Y\smallsetminus\Delta)\to G$, which factors through $\Fie\colon H_1(Y\smallsetminus\Delta,\ZZ)\to G$ since $G$ is abelian. In this paper $Y$ is always $\PP^1$, $\PP^2$ or a nonsingular del Pezzo surface, and $\Delta$ is a configuration of points or lines, hence
\[H_1(Y\smallsetminus\Delta,\ZZ)\cong\frac{\bigoplus_{i=1}^n\ZZ\Span{\Delta_i}}{\Delta},\]
where $\Delta_i$ are the irreducible components of $\Delta$ so that $\Delta=\sum_{i=1}^n\Delta_i$.

Now suppose $P$ is a point of intersection of two or more branch lines $l_1,\dots,l_k$. Then $X$ is nonsingular over $P$ if and only if there are exactly two lines intersecting transversely at $P$, and
\begin{equation}\label{eqn!nonsingular}\Span{\Fie(l_1),\Fie(l_2)}=\Span{\Fie(l_1)}\oplus\Span{\Fie(l_2)}\subset G.\end{equation}
If we blow up $P$ to obtain the exceptional curve $E$, then there is an induced group homomorphism on the blow up, which we still call $\Fie$, defined by
\begin{equation}\label{eqn!exceptional}\Fie(E)=\sum_{i=1}^k\Fie(l_i).\end{equation}
Note that the induced cover on the blow up may be unramified over $E$, and that the blow up may or may not resolve the singularity on $X$ over $P$. 

We have
\[\fie_*\Oh_X=\bigoplus_{\chi\in G^*}\eL^{-1}_\chi,\]
where $G^*=\Hom(G,\CC^*)$ is the group of characters of $G$, and $G$ acts on $\eL_\chi^{-1}$ with character $\chi$. Fix $\chi\in G^*$ and choose a generator for the image of $\chi$ in $\CC^*$, so that we view $\chi\colon G\to\ZZ/d$, where $d$ is the order of $\chi$. Then the following formula determines the eigensheaf $\eL_\chi=\Oh_Y(L_\chi)$:
\begin{equation}\label{eqn!eigensheaf}dL_\chi=\sum_{i=1}^n\left(\chi\circ\Fie(\Delta_i)\right)\Delta_i.\end{equation}
If $X$ is nonsingular, the ramification formula for $K_X$ gives
\begin{equation}\label{eqn!ramification}K_X=\fie^*\left(K_Y+\sum_{i=1}^n\left(1-\frac1d_i\right)\Delta_i\right),\end{equation}
where $d_i$ is the order of $\Fie(\Delta_i)$ in $G$, and
\begin{equation}\label{eqn!pushforward}\fie_*\omega_X=\bigoplus_{\chi\in G^*}\omega_Y(L_\chi).\end{equation}

\subsection{The Kulikov surface}\label{sec!kulikovsurf}
\begin{figure}
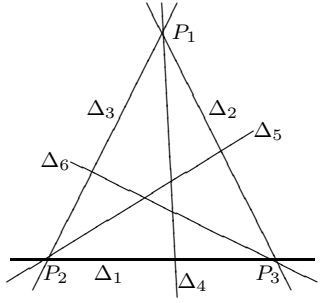

\def\objectstyle{\scriptstyle}
\[
\xy
(0,0)*{}="A"; (30,0)*{}="B"; (15,30)*{}="C";
(15,0)*{}="halfAB"; (20,15)*{}="halfBC"; (10,15)*{}="halfCA";
"A"-(5,0); "B"+(5,0) **\dir{-};
"B"+(2,-4); "C"+(-2,4) **\dir{-};
"C"+(2,4); "A"+(-2,-4) **\dir{-};
"A"-(5.4,3); (27,17) **\dir{-};
"B"+(5.4,-3); (3,13) **\dir{-};
"C"+(0,4); (17,-5) **\dir{-};
"A"+(1,-2)*{P_2};
"B"+(-1,-2)*{P_3};
"C"+(3,0)*{P_1};
"halfAB"+(-7,-2)*{\Del_1};
"halfBC"+(3,5)*{\Del_2};
"halfCA"+(-3,5)*{\Del_3};
"halfAB"+(4,-3)*{\Del_4};
"halfBC"+(9,2)*{\Del_5};
"halfCA"+(-9,-2)*{\Del_6};
\endxy
\]
\caption{The Kulikov configuration}\label{fig!configs}
\end{figure}

In this paragraph we construct the Kulikov surface. We also define the groups $G^0$, $G^1$, $G^2$ which appear throughout this article.
Let $\fie\colon\Xhat\to\PP^2$ be the $(\ZZ/3)^5$-cover of $\PP^2$ branched in the Kulikov configuration $\Delta=\sum_{i=1}^6\Del_i$ of Figure \ref{fig!configs}. Then $\fie$ is determined by the group homomorphism $\Phi\colon H_1(\PP^2\smallsetminus\Delta)\to(\ZZ/3)^6$, represented by the $6\times6$ matrix
\begin{equation}\label{eqn!kulikovFie}
\Fie=\begin{pmatrix}
0&1&0&2&0&0\\
0&0&1&2&0&0\\
\hline
0&0&1&0&2&0\\
1&0&0&0&2&0\\
\hline
1&0&0&0&0&2\\
0&1&0&0&0&2
\end{pmatrix}\end{equation}
of rank $5$. Since $\Delta$ is supported on six lines, $\Fie$ is unique up to
choice of generators for $(\ZZ/3)^6$. Thus $\Xhat$ is maximal, which means any other $(\ZZ/3)^k$-cover of $\PP^2$ branched in $\Delta$ factorises $\fie$.

\begin{dfn}\label{def!G1}\rm Let $G^1\isom(\ZZ/3)^5$ denote the the image of $\Fie$ in $(\ZZ/3)^6$. Then $G^1$ is generated by the columns of $\Fie$,
and we label these $\del_1,\del_2,\del_3,\omega_1,\omega_2,\omega_3$ respectively. Thus $\Xhat$ is a $G^1$-cover of $\PP^2$ branched in the Kulikov config\-uration.
\end{dfn}

The Kulikov surface $\psi\colon X\to\PP^2$ is a $(\ZZ/3)^2$-cover of $\PP^2$ branched in the Kulikov line configuration, and first appeared in \cite{Kul}. The cover is determined by the homomorphism
$\Psi\colon H_1(\PP^2\smallsetminus\Delta)\to(\ZZ/3)^2$, given by the matrix 
\begin{equation}\label{eqn!kulikovPsi}\Psi=\begin{pmatrix}1&1&1&0&1&2\\0&0&0&1&1&1\end{pmatrix}.\end{equation}

Now $\fie\colon\Xhat\to\PP^2$ factors through the Kulikov surface, giving commutative diagrams
\begin{equation}\label{eqn!commdiag}\xymatrix
{\Xhat\ar[r]^\fie\ar[d]_{\psihat}&\PP^2&(\ZZ/3)^6\ar[d]_\Psihat&{H_1(\PP^2\smallsetminus\Delta,\ZZ)}\ar[l]_<<<<\Fie\ar[dl]^\Psi\\
X\ar[ur]_\psi&&(\ZZ/3)^2}
\end{equation}
where \[\Psihat=\begin{pmatrix}0&0&1&1&0&1\\1&1&2&0&0&2\end{pmatrix},\] and $\psihat\colon\Xhat\to X$ is Galois \'etale with group $G^2=\Ker(\Psihat)\cap\Ima(\Fie)$.

\begin{dfn}\label{def!G0G2}\rm
We define $G^2=\Ker(\Psihat)\cap G^1$, so that $G^2=\Span{g_1,g_2,g_3}\cong(\ZZ/3)^3$ where
\begin{equation}\label{eqn!G2kulikov}
\begin{split}g_1&={}^t\begin{pmatrix}0,&0,&1,&0,&1,&2\end{pmatrix},\\
      g_2&={}^t\begin{pmatrix}1,&2,&0,&0,&1,&0\end{pmatrix},\\
      g_3&={}^t\begin{pmatrix}1,&0,&1,&2,&0,&0\end{pmatrix}.\end{split}\end{equation}
We also define $G^0\subset G^1$ to be the subgroup generated by $\xi_i=\del_{i-1}^2\omega_i\omega_{i+1}$ for
$i=1,2,3$. Then $G^0$ is isomorphic to $(\ZZ/3)^2$ because $\xi_1\xi_2\xi_3=1$. We note the relations $g_i=\xi_{i+1}\omega_{i+1}$ for $i=1,2,3$.
\end{dfn}

Note that by \eqref{eqn!nonsingular}, $X$ is singular as a cover of $\PP^2$, so we blow up the three singular points $P_1$, $P_2$, $P_3$ of $\Delta$, and view $X$ as a nonsingular $(\ZZ/3)^2$-cover of the degree $6$ del Pezzo surface $Y$. Then computations using (\ref{eqn!exceptional}--\ref{eqn!pushforward}) show that the Kulikov surface $X$ is a minimal surface of general type with $p_g=0$ and $K^2=6$.

\section{Moduli spaces}\label{sec!moduli}
\subsection{Elliptic curves as abelian covers}\label{sec!curves}
This digression on elliptic curves and abelian covers is required for the Inoue construction of the subsequent section.

Let $\fie\colon E\to\PP^1$ be a Galois $(\ZZ/3)^2$-cover branched in $\Delta=P_0+P_1+P_\infty$. Note that $E$ is the (unique) Fermat cubic curve, see Remark \ref{rmk!fermat}. There are several ways to study $\fie$, but for consistency we prefer to use the techniques of abelian covers. Now $\Fie\colon H_1(\PP^1\smallsetminus\Delta)\to(\ZZ/3)^2$ is unique up to choice of generators, and is given by
\[\Fie=\begin{pmatrix}1&2&0\\0&2&1\end{pmatrix}.\]
Using \eqref{eqn!eigensheaf} to study the various eigensheaves and quotients of the cover, we establish the following diagram:
\begin{equation}\begin{split}\label{eqn!diagcurves}
\xymatrix
{&E\ar[dd]^\fie\ar[dl]_{\pi_1}\ar[dr]^{\pi_2}&\\
E\ar[dr]_{\rho_1}&&\PP^1\ar[dl]^{\rho_2}\\
&\PP^1&}
\end{split}\end{equation}
Here $\pi_1$ is the isogeny induced by taking the quotient by the group $\Span{\eta}$ of translations by the $3$-torsion point $\eta$. There are three possibilities for $\pi_2$ corresponding to the remaining nontrivial subgroups of $(\ZZ/3)^2$. We choose the subgroup $\Span{\omega}$ which fixes the three points $\fie^{-1}(P_1)$ on $E$, and we label these points $0,\eta,2\eta$. Since we have now fixed the origin of $E$, we may consider $\omega=e^{2\pi i/3}$ as a rotation of $E$, under which $\eta$ is fixed. The map $\rho_1$ is the $\ZZ/3$-cover of $\PP^1$ branched in $\Delta$, and $\rho_2$ is the $\ZZ/3$-cover of $\PP^1$ branched over $P_0+P_\infty$. Note that in terms of the original basis for $(\ZZ/3)^2$,
\begin{equation}\label{eqn!etaomega}\eta=\begin{pmatrix}1\\2\end{pmatrix}\text{ and }\omega=\begin{pmatrix}2\\2\end{pmatrix}.\end{equation}
\begin{rmk}\label{rmk!fermat}\rm
Let $E\colon(a^3+b^3+c^3=0)\subset\PP^2$ be the Fermat cubic curve. Then one way to write the $(\ZZ/3)^2$-action is:
\[
    \eta\colon
    \begin{pmatrix}
      a \\ b \\ c
    \end{pmatrix}
    \mapsto
    \begin{pmatrix}
      a \\ \omega b \\ \omega^2c
    \end{pmatrix},\quad
    \omega\colon
    \begin{pmatrix}
      a \\ b \\ c
    \end{pmatrix}
    \mapsto
    \begin{pmatrix}
      \omega a \\ b \\ c
    \end{pmatrix}.
  \]
The map $\pi_2$ is the projection to the line $(a=0)\subset\PP^2$ with coordinates $(b,c)$, so that $\rho_2\colon(b,c)\mapsto(b^3,c^3)$. Then the three ramification points $(0,1,-1),(0,1,-\omega),(0,1,-\omega^2)$ of $\pi_2$ are $0,\eta,2\eta$ respectively. The image of $\pi_1$ is $E\colon(s^3+x^2\xprime+xx'^2=0)\subset\PP^2$, where $x=b^3$, $\xprime=c^3$, $s=abc$. This is the Fermat cubic again, via the relation $s^3=a^3x\xprime$.
\end{rmk}

\subsection{The Inoue-type construction}\label{sec!inoue}
In this paragraph we give an alternative description of the Kulikov surface, which is called the Inoue-type construction, after \cite{Inoue}. This sets the Kulikov surface into context alongside the Burniat--Inoue surfaces.

First note that after blowing up $P_1,P_2,P_3$ we may think of the Kulikov configuration in Figure \ref{fig!configs} as an element of the pencil \[Y_{1,1,1}\subset\PP^1\times\PP^1\times\PP^1\] of degree $6$ del Pezzo surfaces, with equation $\la x_1x_2x_3=\mu\xprime_1\xprime_2\xprime_3$, where $x_i,\xprime_i$ are the coordinates on the $i$th factor. Now let $E\xrightarrow{\pi_2}\PP^1\xrightarrow{\rho_2}\PP^1$ be a $(\ZZ/3)^2$-cover of $\PP^1$ as described in Section \ref{sec!curves}. Taking a direct product of three copies of $E$, we obtain the following diagram:
\begin{equation}\begin{split}\label{eqn!inouediag}
\xymatrix@R=1pt
{E\times E\times E\ar[r]^\pi&\PP^1\times\PP^1\times\PP^1\ar[r]^\rho&\PP^1\times\PP^1\times\PP^1\\
\bigcup&\bigcup&\bigcup\\
\Xhat_{\phantom{}}\ar[r]\ar[dddr]_\psihat&Z\cup\bigcup_iZ_i\ar[r]&Y_{1,1,1}\\
&&\\
&&\\
&X\ar[uuur]_\psi&}
\end{split}\end{equation}
where $\pi=\prod\pi_2^i$, $\rho=\prod\rho_2^i$.

Fix an element $Y$ of our pencil of del Pezzo surfaces. Then $\rho^{-1}(Y)$ splits into three components $Z_i$, $i=0,1,2$, each of which has stabiliser isomorphic to $(\ZZ/3)^2$. We fix $Z=Z_0$, with stabiliser $G^0\cong(\ZZ/3)^2$. Now define $\Xhat=\pi^{-1}(Z)$, so that $\Xhat$ is a hypersurface in $E\times E\times E$ of tridegree $(3,3,3)$, with stabiliser $G^1=G^0\oplus(\ZZ/3)^3$. The quotient of $\Xhat$ by $G^1$ is $Y$, and by construction, the branch locus of $\Xhat\to Y$ is exactly the blown up Kulikov configuration. Thus $\Xhat$ is the maximal $(\ZZ/3)^5$-cover of Section \ref{sec!kulikovsurf}.
\begin{lemma}\label{lemma!Inoue} There is a subgroup $G^2\isom(\ZZ/3)^3$ of $G^1$ such that $G^2$ acts freely on $\Xhat$, and the quotient $X$ is a Kulikov surface.
\end{lemma}
\begin{pf} Since $\Xhat\to Y$ is the maximal $(\ZZ/3)^5$-cover of $Y$, we see that the quotient map $(\rho\circ\pi)|_{\Xhat}=\fie$, where $\fie$ is determined by the matrix \eqref{eqn!kulikovFie} of Section \ref{sec!kulikovsurf}. Hence the group $G^1$ is generated by the columns $\del_i$, $\omega_i$ of \eqref{eqn!kulikovFie} as stated in Definition \ref{def!G1}. By Definition \ref{def!G0G2}, the subgroup $G^2$ is generated by
\[g_i=\xi_{i+1}\omega_{i+1}=\del_i^2\omega_{i+1}^2\omega_{i+2},\]
for $i=1,2,3$.

For completeness, we reconcile the above construction of $G^0$ with Definition \ref{def!G0G2}. Observe that we subdivided the rows of matrix \eqref{eqn!kulikovFie} into pairs, each of which corresponds to the action of $G^1$ restricted to the respective factor of $E\times E\times E$. For example, by Definition \ref{def!G0G2},
\[\xi_1={}^t\begin{pmatrix}2,&1,&1,&2,&0,&0\end{pmatrix},\]
and comparing this with the definition of $\eta=\binom12$ from \eqref{eqn!etaomega}, we see that
\[\xi_1=\eta_1^2\eta_2,\]
where $\eta_i$ now denotes translation by the $3$-torsion point $\eta$ on the $i$th factor of the product. We also have
\[\xi_2=\eta_2^2\eta_3,\ \xi_3=\eta_1\eta_3^2.\]
This is called the Inoue-type construction of Kulikov surfaces, cf.~\cite{BC}, \cite{Inoue} for the original Inoue surfaces.$\qed$
\end{pf}

\subsection{The moduli space of Kulikov surfaces}
Using the Inoue-type construction, it is clear that we obtain a pencil of Kulikov surfaces
as the pullback of the pencil of del Pezzo surfaces. In fact we have:
\begin{thm} The pencil of Kulikov surfaces forms a $1$-dimensional ir\-reducible component
of the moduli space of surfaces of general type with $p_g=0$ and $K^2=6$.
\end{thm}
\begin{pf}
Let $X$ be a Kulikov surface. We calculate $h^1(T_X)=h^1(\psi_*T_X)$ using the splitting of
$\psi_*T_X$ into eigensheaves according to Proposition 4.1 of \cite{Par}:
\begin{equation}
\begin{split}\label{eqn!pardinitgt}
(\psi_*T_X)^{\inv}&=T_Y(-\log\Del)\\
(\psi_*T_X)^{(\chi)}&=T_Y(-\log\Del_g:g\in S_\chi)\otimes\eL_\chi^{-1}
\end{split}
\end{equation}
where $S_\chi = \{g\in G | \chi(g)\ne m - 1\}$, $m$ is the order of $g$, and
$\Del_g$ is the sum of the components $\Del_i$ of $\Del$ such that $\Fie(\Del_i)=g$.
For example, when $g=\binom10$ we have
\[\Del_{\binom{1}{0}} = \Del_1 + \Del_2 + \Del_3 = 3H - 2\sum_{i=1}^3 E_i,\]
where $E_i$ is the exceptional curve over $P_i$. The first equality follows from the definition
of $\Psi$ in \eqref{eqn!kulikovPsi} and by applying equation \eqref{eqn!exceptional} to each exceptional curve. The remaining cases are
  \begin{align*}  
    \Del_{\binom{1}{1}} &= \Del_5 + E_3 = H - E_2 + E_3, \\
    \Del_{\binom{2}{1}} &= \Del_6 + E_1 = H - E_3 + E_1, \\
    \Del_{\binom{0}{1}} &= \Del_4 + E_2 = H - E_1 + E_2,
  \end{align*}
and all other $\Del_g$ are zero.

Given the decomposition \eqref{eqn!pardinitgt}, it is clear that
\begin{equation}\label{eqn!h1decomp}h^i(T_X)=\sum_\chi h^i((\psi_*T_X)^{(\chi)}).\end{equation}
We use several methods, which we describe subsequently, to calculate the value of
 $h^2((\psi_*T_X)^{(\chi)})$ for each $\chi$. For convenience, we present the results in Table \ref{tab!kulmod}.
Now since $h^0(T_X)=0$ for a surface of general type, we have
\[\chi(T_X)=-h^1(T_X)+h^2(T_X)=2K_X^2-10\chi(\Oh_X)=2.\]
So using this equality, Table \ref{tab!kulmod} and equation \eqref{eqn!h1decomp} we get $h^1(T_X)=1$, which proves the theorem.
$\qed$
\end{pf}
\begin{table}[htbp]
\small
\begin{equation*}
  \renewcommand{\arraystretch}{1.4}
  \begin{array}{|c|c|c|c|c|c|}
    \hline
     & \chi & A=K_Y+L_\chi & \{\Del_g:g\in S_\chi\} & h^2((\psi_*T_X)^{(\chi)})\\
    \hline
    (0) & (0,0) & -3H+\sum_i E_i & \Del & 0\\
    \hline
        & (1,0) & -H + E_1 & \Del_1,\Del_2,\Del_3,\Del_4,\Del_5,E_2,E_3 &\\
    (a) & (1,2) & -H + E_2 & \Del_1,\Del_2,\Del_3,\Del_5,\Del_6,E_1,E_3 & 1\\
        & (1,1) & -H + E_3 & \Del_1,\Del_2,\Del_3,\Del_4,\Del_6,E_1,E_2 &\\
    \hline
    (b)	& (0,1) & -2H + \sum_i E_i & \Del_1,\dots,\Del_6,E_1,E_2,E_3 & 0\\
    \hline
    (c)	& (0,2) & -H + \sum_i E_i & \Del_1,\Del_2,\Del_3 & 0\\
    \hline
    	& (2,2) & -E_1 & \Del_5,\Del_6,E_1,E_3 &\\
    (d)	& (2,0) & -E_2 & \Del_4,\Del_6,E_1,E_2 & 0\\
    	& (2,1) & -E_3 & \Del_4,\Del_5,E_2,E_3 &\\
    \hline
  \end{array}
\end{equation*}
\caption{Decomposition of $\psi_*T_X$.}\label{tab!kulmod}
\end{table}

The remainder of this section deals with Table \ref{tab!kulmod}. The table
contains the data required to write down the eigenspace decomposition of $\psi_*T_X$
according to \eqref{eqn!pardinitgt}. We also remark that the divisors $\Del_g$ are always
composed of a number of disjoint rational curves. Thus there is no need to distinguish between
logarithmic poles along $D_1$ and $D_2$, and logarithmic poles
along $D_1+D_2$ for any curves $D_1,D_2$. The data are grouped into cases $(0),(a),\dots,(d)$ by symmetry
considerations, avoiding repetition in our calculations.

First observe that by Serre duality,
\begin{equation}\label{eqn!SD}h^i((\psi_*T_X)^{(\chi)})=h^{2-i}(\Omega_Y(\log\Del_g:g\in S_\chi)(A)),\end{equation}
where $A=K_Y+L_\chi$. Thus the foundation for our calculations of the various $h^1((\psi_*T_X)^{(\chi)})$
is the standard logarithmic residue sequence twisted by $A$,
\begin{equation}\label{eqn!resseq}
0\to\Omega_Y(A)\to\Omega_Y(\log\Del_g\colon g\in S_\chi)(A)\to\bigoplus_{g\in S_\chi}\Oh_{\Del_g}(A)\to0,
\end{equation}
and its associated long exact sequence of cohomology, with connecting homo\-morphism
\[\del\colon \bigoplus_{g\in S_\chi}H^0(\Oh_{\Del_g}(A))\to H^1(\Omega_Y(A)).\]
Now if $A=0$, then the image of $\del$ is $\Span{c_1(\Del_g)\colon g\in S_\chi}$. Thus it is possible to
calculate the rank of $\del$ directly, because of the Lefschetz $(1,1)$ Theorem. In turn this gives us the
value of $h^0(\Omega_Y(\log\Del_g\colon g\in S_\chi)(A))$ using the long exact sequence.
The crucial part of our proof is to adapt this argument to situations where $A$ is nontrivial.

Note that it is possible to calculate the value of $h^1((\psi_*T_X)^{(\chi)}$ in each case, using the Riemann--Roch
theorem together with the logarithmic residue sequence \eqref{eqn!resseq}. However, this is not needed for
the proof so we omit this calculation.

\subsubsection{The basic method}
In this paragraph we compute the last column of Table \ref{tab!kulmod} in cases $(0)$, $(b)$ and $(d)$.
We first give the set up: suppose that $|{-}A|$ has no fixed part, so there is a morphism $\Oh_Y(A)\to\Oh_Y$ which is nonzero on each irreducible component of the branch divisor $\Delta$. Then we obtain the following commutative diagram:
\[\xymatrix
{0\ar[r]&\Omega_Y(A)\ar[r]\ar[d]&\Omega_Y(\log\Del_g\colon g\in S_\chi)(A)\ar[r]\ar[d]&\bigoplus_{g\in S_\chi}\Oh_{\Del_g}(A)\ar[d]\ar[r]&0\\
0\ar[r]&\Omega_Y\ar[r]&\Omega_Y(\log\Del_g\colon g\in S_\chi)\ar[r]&\bigoplus_{g\in S_\chi}\Oh_{\Del_g}\ar[r]&0
}\]
The connecting homomorphisms for the corresponding long exact sequences fit into the commutative square
\[\xymatrix
{\bigoplus_{g\in S_\chi}H^0(\Oh_{\Del_g}(A))\ar[d]\ar[r]^>>>>>{\del}\ar[dr]^\al&H^1(\Omega_Y(A))\ar[d]^{\be}\\
\bigoplus_{g\in S_\chi}\CC_{\Del_g}\ar[r]^{c_1}&H^1(\Omega_Y)
}\]
where $\CC_{\Del_g}$ denotes $H^0(\Oh_{\Del_g})$, $c_1$ is the Chern class map, $\del$ is the connecting homomorphism whose rank we wish to calculate, and $\be$ has kernel $H^0(\Omega_Y|_D)$,
where $D\in|{-}A|$.

As an illustration, consider case $(b)$, so that $\chi=(0,1)$ and $A=-2H+E_1+E_2+E_3$. Clearly $A$ has no fixed part, and moreover 
\[\bigoplus_{g\in S_{(0,1)}}H^0(\Oh_{\Del_g}(A))=\CC_{\Del_1}\oplus\CC_{\Del_2}\oplus\CC_{\Del_3}.\]
Now since $\Del_1$, $\Del_2$ and $\Del_3$ are linearly independent in $\Pic Y$, $\rank\al=3$. Thus $\rank\del=3$, and so $h^0\left(\Omega_Y(\log\Del_g\colon g\in S_{(0,1)})(A)\right)=0$.

Case $(0)$ is also treated using this method, since $\bigoplus_{g\in G}H^0(\Oh_{\Del_g}(K_Y))=0$, giving $h^0\left(\Omega_Y(\log\Del)(K_Y)\right)=0$.
For case $(d)$, suppose $\chi=(2,2)$, so that $A={-}E_1$. Then we have a natural inclusion 
\[H^0(\Omega_Y(\log\Del_g\colon g\in S_{(2,2)})(-E_1))\subseteq H^0(\Omega_Y(\log\Del_g\colon g\in S_{(2,2)})),\]
and we can show the right hand side vanishes by using the standard log\-arithmic residue sequence.

\subsubsection{The contraction lemma and fibration method}
This paragraph deals with case $(a)$, where the basic method does not work. We contract some branch divisors on $Y$ and then calculate directly.
A more general version of the following lemma appears in \cite{BCII}, \cite{BCIII}. We quote a simplified version, which is sufficient for our purposes:
\begin{lemma}\label{lem!cont}
Let $\pi\colon Z\to\CC^2$ be the blowup of $\CC^2$ at the origin $O$, with exceptional curve $E$. Suppose $L_1$, $L_2$ are lines through $O$ with distinct tangents at $O$, and let $D_i$ be the strict transform of $L_i$. Then
\[\pi_*\Omega_Z(\log D_1, \log D_2, \log E)=\Omega_{\CC^2}(\log L_1,\log L_2).\]
\end{lemma}

In case $(a)$, with $\chi=(1,0)$, we have $A=-H+E_1$. Note that the basic method does not work here: the rank of $\al$ is $3$, which is not maximal, so we do not have sufficient information to calculate $\rank\del$. Instead we contract $\Del_3$ and $E_3$ via $\pi\colon Y\to Q=\PP^1\times\PP^1$, and let $f\colon Q\to\PP^1$
be the projection to the second factor. Then by Lemma \ref{lem!cont},
\[\pi_*(\Omega_Y(\log\Del_g\colon g\in S_{(1,0)})(A))=\Omega_Q(\log(F_1+F_2+F_3),\log(B_1+B_2))(-F_2),\]
where $F_1=\pi_*\Del_2$, $F_2=\pi_*A$, $F_3=\pi_*E_2$ are fibres of $f$, and $B_1=\pi_*\Del_1$, $B_2=\pi_*\Del_5$ are sections.
Now, the standard short exact sequence
\[0\to f^*\omega_{\PP^1}\to\Omega_Q\to\omega_{Q/\PP^1}\to0\]
specialises to
\[0\to f^*\omega_{\PP^1}(F)\to\Omega_Q(\log F, \log B)\to\omega_{Q/\PP^1}(B)\to0.\]
Thus twisting by $-F_2$ and observing that $f^*\omega_{\PP^1}=\Oh_Q(-2F)$ and $\omega_{Q/\PP^1}=\Oh_Q(-2B)$, we get
\[h^0(\Omega_Y(\log\Del_g\colon g\in S_{(1,0)})(A))=h^0(\Oh_Q)=1.\]

\subsubsection{Final case (c)}
In case $(c)$, $A=-H+\sum E_i$, so $|{-}A|$ is empty and it is not possible to use the basic method. We can not use the contraction lemma either, because when we contract $\Del_i$, we still have poles along $E_i$, and further contractions are complicated. Instead we use the natural inclusion
\[H^0(\Omega_Y(\log(\Del_1,\Del_2,\Del_3)(A))\subseteq H^0(\Omega_Y(B)),\]
where $B=A+\Del_1+\Del_2+\Del_3$. Then the exact sequence
\[0\to\si^*\Omega_{\PP^2}\to\Omega_Y\to\bigoplus_{i=1}^3\Oh_{E_i}(-2)\to0\]
gives $H^0(\Omega_Y(B))=H^0(\si^*\Omega_{\PP^2}(B))$, and the pullback of the dual Euler sequence to $Y$
\[0\to\si^*\Omega_{\PP^2}\to\Oh_Y(-H)^3\to\Oh_Y\to0\]
gives $H^0(\si^*\Omega_{\PP^2}(B))=H^0(\Oh_Y(-A)^3)=0$. Hence 
\[h^0(\Omega_Y(\log\Del_g\colon g\in S_{(0,2)})(A))=0.\]

\subsection{Degenerate Kulikov surfaces}\label{sec!degenerate}
In this paragraph we constuct three degenerate Kulikov surfaces arising as $(\ZZ/3)^2$-covers of special elements of the pencil of del Pezzo surfaces
\[Y_{1,1,1}\colon(\la x_1x_2x_3=\mu x'_1x'_2x'_3)\subset\PP^1\times\PP^1\times\PP^1.\]
These degenerate surfaces are members of the boundary of the compactified moduli space, in the sense of \cite{KSB}.
\begin{figure}[htbp]
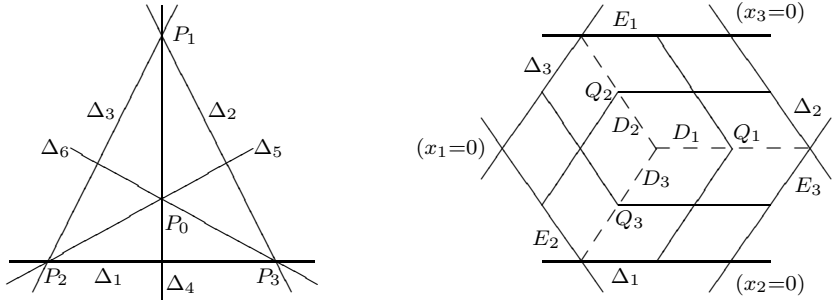

\def\objectstyle{\scriptstyle}
\[
\xy
(0,0)="O";
"O"+(0,0)*{}="OA"; "O"+(30,0)*{}="OB"; "O"+(15,30)*{}="OC";
"O"+(15,0)*{}="OhalfAB"; "O"+(20,15)*{}="OhalfBC"; "O"+(10,15)*{}="OhalfCA";
"OA"-(5,0); "OB"+(5,0) **\dir{-};
"OB"+(2,-4); "OC"+(-2,4) **\dir{-};
"OC"+(2,4); "OA"+(-2,-4) **\dir{-};
"OA"-(5.4,3); "O"+(27,15) **\dir{-};
"OB"+(5.4,-3); "O"+(3,15) **\dir{-};
"OC"+(0,4); "O"+(15,-5) **\dir{-};
"O"+(17,5)*{P_0};
"OA"+(1,-2)*{P_2};
"OB"+(-1,-2)*{P_3};
"OC"+(3,0)*{P_1};
"OhalfAB"+(-7,-2)*{\Del_1};
"OhalfBC"+(3,5)*{\Del_2};
"OhalfCA"+(-3,5)*{\Del_3};
"OhalfAB"+(2.5,-3)*{\Del_4};
"OhalfBC"+(9,0)*{\Del_5};
"OhalfCA"+(-9,0)*{\Del_6};
(60,0)="Op";
"Op"+(10,0)*{}="OpA"; "Op"+(30,0)*{}="OpB";
"Op"+(40,15)*{}="OpC"; "Op"+(30,30)*{}="OpD";
"Op"+(10,30)*{}="OpE"; "Op"+(0,15)*{}="OpF";
"Op"+(20,0)*{}="OphalfAB"; "Op"+(35,7.5)*{}="OphalfBC";
"Op"+(35,22.5)*{}="OphalfCD"; "Op"+(20,30)*{}="OphalfDE";
"Op"+(5,22.5)*{}="OphalfEF"; "Op"+(5,7.5)*{}="OphalfAF";
"OpA"-(5,0); "OpB"+(5,0) **\dir{-};
"OpB"+(-3,-4); "OpC"+(3,4) **\dir{-};
"OpC"+(3,-4); "OpD"+(-3,4) **\dir{-};
"OpE"-(5,0); "OpD"+(5,0) **\dir{-};
"OpF"+(-3,-4); "OpE"+(3,4) **\dir{-};
"OpA"+(3,-4); "OpF"+(-3,4) **\dir{-};
"Op"+(20,15)="Centre";
"OpA"; "Centre" **\dir{--};
"OpC"; "Centre" **\dir{--};
"OpE"; "Centre" **\dir{--};
"OphalfAB"; "Op"+(30,15) **\dir{-};
"Op"+(30,15); "OphalfDE"; **\dir{-};
"OphalfBC"; "Op"+(15,7.5) **\dir{-};
"Op"+(15,7.5); "OphalfEF"; **\dir{-};
"OphalfCD"; "Op"+(15,22.5) **\dir{-};
"Op"+(15,22.5); "OphalfAF"; **\dir{-};
"OpA"+(6,-2)*{\Del_1};
"OpC"+(0,-5)*{E_3};
"OpC"+(0,5.5)*{\Del_2};
"OpE"+(6,2)*{E_1};
"OpE"+(-5.5,-4)*{\Del_3};
"OpA"+(-4.5,3)*{E_2};
"OpB"+(5,-3)*{(x_2=0)};
"OpD"+(5,3)*{(x_3=0)};
"OpF"+(-7,0)*{(x_1=0)};
"Centre"+(4,2)*{D_1};
"Centre"+(-4,3)*{D_2};
"Centre"+(0,-4)*{D_3};
"Op"+(32,17)*{Q_1};
"Op"+(12.5,22.5)*{Q_2};
"Op"+(16.5,5.5)*{Q_3};
\endxy
\]
\caption{The complete quadrangle and the cube}\label{fig!compquad}
\end{figure}
\subsubsection{A hyperelliptic surface}
When $\la=\mu=1$, the three lines $\Del_4$, $\Del_5$ and $\Del_6$ of the Kulikov configuration meet in a point $P_0$. This configuration is called the complete quadrangle and it is illustrated on the left side of Figure \ref{fig!compquad}. The $(\ZZ/3)^2$-cover $X\to Y$ has an elliptic singularity of degree $9$ over $P_0$. If we blow up $Y$ at $P_0$, then the induced $(\ZZ/3)^2$-cover $\Xtilde$ is the resolution.

Note that the pencil of conics in $\PP^2$ passing through $P_0,\dots,P_3$ gives rise to a pencil of Fermat cubics on $\Xtilde$, so $\Xtilde$ is an elliptic surface. Moreover, the minimal model of $\Xtilde$ is the $(\ZZ/3)^2$-cover of $\PP^1\times\PP^1$ obtained by contracting the strict transforms of $\Del_4$, $\Del_5$ and $\Del_6$. The branch locus is transformed to six lines (three in each ruling), and the minimal model of $\Xtilde$ is the hyperelliptic surface $(E\times E)/(\ZZ/3)^2$.

One can also see this degeneration in terms of the Inoue-type construction of Section \ref{sec!inoue}. The action of $G^2$ is free on $E\times E\times E$ outside the orbit of the point $(0,0,0)$, which
has stabiliser $\ZZ/3$ generated by $g_1g_2g_3$. The deg\-enerate Kulikov surface described above is obtained when the hypersurface $\Xhat_{3,3,3}\subset E\times E\times E$ contains the point $(0,0,0)$.

\subsubsection{Two reducible surfaces}
When $\mu=0$ or $\la=0$, the del Pezzo surface $Y$ breaks into three copies of $\PP^1\times\PP^1$. We only consider the case $\mu=0$ as shown on the right side of Figure \ref{fig!compquad}, because $\la=0$ has a similar treatment by symmetry. The components of $Y$ are called $Y_i\colon(x_i=0)$, and $Y$ has normal crossing singularities along the lines $D_k\colon(x_i=x_j=0)$ for $i,j,k=\{1,2,3\}$, indicated by dotted lines in the figure. The divisors $\Del_4$, $\Del_5$ and $\Del_6$ break into pairs of lines, for example $\Del_4$ consists of the two lines joining $E_1$ to $\Del_1$, intersecting in a single point $Q_1$. For clarity, we do not label $\Del_4$, $\Del_5$, $\Del_6$ in the figure.

We construct each $(\ZZ/3)^2$-cover $\psi_i\colon X_i\to Y_i$ separately, before gluing them together to give the degenerate Kulikov surface. For $i=1$, we see from the figure that $\psi_1$ is branched over the four lines $E_1$, $\Del_2$, $\Del_4$ and $\Del_5$. Recall from Section \ref{sec!surfaces} that the $(\ZZ/3)^2$-cover of $Y_1$ branched in $\Del$ is governed by the group homomorphism $\Psi_1\colon H_1(Y_1\smallsetminus\Del,\ZZ)\to(\ZZ/3)^2$. Thus by \eqref{eqn!kulikovPsi}, $\psi_1$ must also be branched over $D_1$ and $D_2$, with $\Psi_1(D_1)=-\Psi_1(E_1)-\Psi_1(\Del_5)=\binom01$, $\Psi_1(D_2)=-\Psi_1(\Del_2)-\Psi_1(\Del_4)=\binom22$. Then by equations \eqref{eqn!nonsingular} and \eqref{eqn!exceptional}, $X_1$ has three singularities of type $\frac13(1,1)$ over $Q_2$, and three of type $\frac13(1,2)$ over $Q_1$. If we blow up $Q_2$ and contract the strict transforms of $D_2$ and $\Del_5$, the induced cover ${\Xtilde}_1\to\PP^2$ is a rational surface with three $\frac13(1,2)$ singularities.

The other components of $X$ have similar singularities over the points $Q_i$. Perhaps the most interesting aspect of this example is that when we glue the components of $X$ back together, we see that over each $Q_i$, we must attach a $\frac13(1,2)$ point to a $\frac13(1,1)$ point. This gives rise to orbifold normal crossing singularities, which can be expressed locally as
\[(xy=0)\subset\frac13(1,2,1),\]
where the double curve is given by $(x=y=0)$.

\section{Fundamental groups and homotopy equivalence}\label{sec!homotopy}
In this section we give an explicit description of the fundamental group of the Kulikov surface. Using
this description we show that the $1$-dimensional irreducible component of the moduli space constructed
in the previous section is actually closed under homotopy equivalence. Hence the Kulikov surfaces form a
connected component of the moduli space.

\subsection{Fundamental group and homology}
Using the Inoue-type construction of Section \ref{sec!inoue}, it is
possible to lift the action of $G^2=\Span{g_1,g_2,g_3}$ to the Fermat cubic curves.
Indeed, let $z_i$ be a uniformising parameter for the $i$th curve in the direct product
\eqref{eqn!inouediag}. Recall that in the proof of Lemma \ref{lemma!Inoue}, we showed that $\xi_i=\eta_i^2\eta_{i+1}$
for $i=1,2,3$. Now by Definition \ref{def!G0G2}, we have $g_i=\xi_{i+1}\omega_{i+1}$. Hence 
the action of $G^2$ lifts to the product $E\times E\times E$ as
\begin{equation*}
  g_1(z)=\begin{pmatrix}z_1\\\omega z_2+2\eta_2\\z_3+\eta_3\end{pmatrix},\
  g_2(z)=\begin{pmatrix}z_1+\eta_1\\z_2\\\omega z_3+2\eta_3\end{pmatrix},\
  g_3(z)=\begin{pmatrix}\omega z_1+2\eta_1\\z_2+\eta_2\\z_3\end{pmatrix}.
\end{equation*}

Now write $e_i,e'_i$ for the lattice generators of the $i$th factor in $E\times E\times E$.
In particular, we fix $e_i=1$, $e'_i=\omega$ so that the $3$-torsion points $\eta_i=\frac13(2e_i+e'_i)$
and $2\eta_i=\frac13(e_i+2e'_i)$ are fixed under rotation by $\omega$.
We lift $g_i$ to affine transformations $\ga_i$ in $\AA(3,\CC)$, and define $\Gam\subset\AA(3,\CC)$ to be the
subgroup generated by $\ga_i,t_i,t'_i$ for $i=1,2,3$, where $t_i,t'_i$ are affine translations by $e_i,e'_i$
respectively. Note that since $\ga_i^3 = t_{3\eta_{i+2}} = t_{i+2}^2t'_{i+2}$,
where $i$ is taken modulo $3$, $\Gam$ is in fact generated by the $\ga_i$ and $t_i$.
\begin{thm}\label{thm!kulthreehomology}
A Kulikov surface $X$ has $\pi_1(X) = \Gam$ and $H_1(X,\ZZ) = \paren{\ZZ/3}^3$.
\end{thm}
\begin{rmk}\rm Kulikov proved in \cite{Kul} that $\pi_1(X)$ is infinite and not abelian, and that $\Tors X\supset(\ZZ/3)^3$.\end{rmk}
\begin{pf}
By the Inoue-type construction of Section \ref{sec!inoue}, $X$ has an \'etale
 $\paren{\ZZ/3}^3$-cover $\Xhat$, which is a smooth hypersurface of tridegree $(3,3,3)$ in $E\times E\times E$.
  Therefore, by Lefschetz's theorem, $\pi_1(\Xhat) = \pi_1(E\times
  E\times E)$, which is generated by translations $t_i$ and $t'_i$ for $i=1,2,3$.

Now, $\Gam$ acts on the universal cover $\CC^3$ of $E\times E\times E$, and it acts freely on the
universal cover $\Xtilde\subset\CC^3$ of $\Xhat$. By construction $X = {\Xtilde}/\Gam$, so $\Xtilde$
is the universal cover of $X$ and $\pi_1(X) = \Gam$. Clearly, we have also shown that there is an exact sequence
\[1\to\ZZ^6\to\pi_1(X)\to G^2\to1.\]

The group $H_1(X,\ZZ)$ is the abelianisation of $\Gam$, so we calculate its centre $[\Gam,\Gam]$.
Observe that
  \begin{equation}\label{eqn!comrel3-1}
    \begin{split}
      \ga_1\ga_2 &= t_3\ga_2\ga_1 = \ga_2\ga_1(t_3t'_3)^{-1},\\
      \ga_2\ga_3 &= t_1\ga_3\ga_2 = \ga_3\ga_2(t_1t'_1)^{-1},\\
      \ga_3\ga_1 &= t_2\ga_1\ga_3 = \ga_1\ga_3(t_2t'_2)^{-1},
    \end{split}
  \end{equation}
hence all $t_i,t'_i$ are in $[\Gam,\Gam]$. Moreover, $\Gam/{\Span{t_{i},t'_{i}}}$ is abelian, so
$[\Gam,\Gam]=\Span{t_i,t'_i}$ and
\begin{equation*}
  H_1(X,\ZZ) = \Gam^{\text{ab}}\isom\paren{\ZZ/3}^3,
\end{equation*}
where the generators are the residue classes of $\ga_i$ modulo $[\Gam,\Gam]$.$\qed$
\end{pf}

\subsection{The moduli space is closed under homotopy equivalence}
In this paragraph we prove that the moduli space of Kulikov surfaces is actually a connected component of the moduli space of surfaces of general type with $K^2=6$. In fact we prove a stronger result:
\begin{thm}\label{thm:homot-equiv-Kul3}
Any compact surface $S$ which is homotopically equivalent to a Kulikov surface $X$ is itself a Kulikov surface.
\end{thm}

We prove the theorem in several steps, mostly by an explicit analysis of the fundamental group. A similar theorem for Burniat--Inoue surfaces is proved in \cite{BC}. For the remainder of this section, $X$ is a Kulikov surface and $S$ is a compact surface homotopically equivalent to $X$.

Now consider the three subgroups of $\Gam$ defined by
  \begin{align*}
    \Si_1 &= \Span{t_1, \ga_1, t_2, \ga_2, t_3, \ga_3^3},\\
    \Si_2 &= \Span{t_1, \ga_1^3, t_2, \ga_2, t_3, \ga_3},\\
    \Si_3 &= \Span{t_1, \ga_1, t_2, \ga_2^3, t_3, \ga_3}.
  \end{align*}
These are of index $3$ in $\Gam$ and using the commutation relations
(\ref{eqn!comrel3-1}) as well as
\begin{equation}\label{eqn!comrel3-2}
  \ga_it_j=
    \left\{\begin{array}{ll}
    (t_j\ga_i)t_j^2t'_j = t'_jt_j^{-1}(t_j\ga_i)&\text{if }j=i+1,\\
    t_j\ga_i&\text{otherwise,}
    \end{array}\right.
\end{equation}
we see that they are normal subgroups. Thus each $\Si_i$ corresponds to an
\'etale $\ZZ/3$-cover $S_i\to S$.

\begin{lemma} The Albanese variety $\Alb(S_i)$ of each $S_i$ is the Fermat
  cubic curve. In other words, $\Alb(S_i)$ is the unique elliptic curve which
  admits complex multiplication by a primitive cube root of unity.
\end{lemma}
\begin{pf}
We prove the lemma for $S_1$, since the other cases are the same by symmetry. Using commutation
relations (\ref{eqn!comrel3-1}) and (\ref{eqn!comrel3-2}), we have
\[[\Si_1,\Si_1] \supset \Span{t_2^2t'_2, t'_2t_2^{-1}, t_3,t_3'},\]
and since $\Si_1/{\Span{t_2^2t'_2, t'_2t_2^{-1}, t_3,t_3'}}$ is abelian, the reverse inclusion holds. Hence
\[H_1(S_1,\ZZ) = \Si_1^\text{ab}\isom \Span{\coset t_1,\coset\ga_2} \oplus (\ZZ/3)^2,\]
where $\coset t_1,\coset \ga_2$ denote the classes of $t_1,\ga_2$ respectively,
and the torsion summand $(\ZZ/3)^2$ is generated by the classes of
$\ga_1$ and $t_2$. This shows that $\Alb(S_1)$ is an elliptic curve.

Now the $\Gam/\Si_1$-action associated to the cover $S_1\to S$ descends to a nontrivial
$\ZZ/3$-action on $H^0(S_1,\Omega_{S_1})$. Indeed,
\[H^0(S_1,\Omega_{S_1})^{\Gam / \Si_1} = H^0(S,\Omega_S)=H^1(S,\Oh_S)=0,\]
by Hodge theory. Hence $\Alb(S_1)$ admits an automorphism of order $3$.
This proves that $\Alb(S_1)$ is the Fermat cubic curve.$\qed$
\end{pf}

Now write $E_i'=\Alb(S_i)$, and let $\La_i'=\Span{t_i,\ga_{i+1}}$ be the lift of the torsion free
part of $H_1(S_i,\ZZ)$ to $\Gam$. Then $\Span{t_i, \ga_{i+1}^3} = \Span{t_i, t_i'}$ is an index
$3$ normal subgroup of $\La_i'$, which we call $\La_i$. Viewed as a lattice, $\La_i$
corresponds to an elliptic curve $E_i$ and a $\ZZ/3$-cover $E_i\to E_i'$. Since $E_i$ and $E_i'$
are isogenous, $E_i$ is also a Fermat cubic curve.

\begin{prop}\label{prop!albanese} Let $\Shat\to S$ be the maximal abelian cover of $S$. Then
\begin{enumerate}
\item[(1)] the Albanese variety of $\Shat$ is the product of elliptic curves $E_1\times E_2\times E_3$;
\item[(2)] the Albanese map $\alb\colon\Shat\to E_1\times E_2\times E_3$ is a birational morphism onto its image
 $A_{\Shat}$, and $A_{\Shat}$ is a hypersurface of tridegree $(3,3,3)$.
\end{enumerate}
\end{prop}

\begin{pf} Define $\La=\La_1\oplus\La_2\oplus\La_3$, and let $\psihat\colon\Shat\to S$ be
the \'etale $(\ZZ/3)^3$-cover of $S$ with fundamental group $\La\lhd\Gam$. Note that $H_1(S,\ZZ)=(\ZZ/3)^3$, so
$\Shat\to S$ is maximal and thus $\psihat$ factors through each $S_i$. Now composing these factorisations with
$\alb\colon S_i\to E_i'$ and taking the direct product, we obtain a morphism
  \begin{equation*}
    f \colon \Shat \to E_1'\times E_2'\times E_3'.
  \end{equation*}
Moreover, $f$ factors through $E_1\times E_2\times E_3$ via the $\ZZ/3$-covers $E_i\to E_i'$.
Since $\pi_1(E_1\times E_2\times E_3)=\pi_1(\Shat)$, this factorisation of $f$ is the Albanese map of $\Shat$.
This proves part $(1)$ of the proposition.

To prove part $(2)$ we compare $H^*(\Shat,\ZZ)$ with $H^*(\Xhat,\ZZ)$, where $\Xhat$ is the maximal $(\ZZ/3)^3$-cover of a Kulikov surface $X$. First let us factor the Albanese map of $\Shat$ through its image $A_{\Shat}$ as follows:
\begin{equation}\label{eqn!factorShat}\Shat\xrightarrow{a_{\Shat}}A_{\Shat}\xrightarrow{i_{\Shat}}\Alb(\Shat)=E_1\times E_2\times E_3.\end{equation}
Then recall from Section \ref{sec!inoue} that $\Xhat$ is a hypersurface of tridegree $(3,3,3)$ in $E\times E\times E$. Now following the procedure described in part $(1)$, we may also factor the Albanese map of $\Xhat$ through its image $A_{\Xhat}$:
\begin{equation}\label{eqn!factorXhat}\Xhat\xrightarrow{a_{\Xhat}}A_{\Xhat}\xrightarrow{i_{\Xhat}}\Alb(\Xhat)=E\times E\times E.\end{equation}
This time $a_{\Xhat}$ is actually an isomorphism.

Now the factorisations \eqref{eqn!factorShat} and \eqref{eqn!factorXhat} are both constructed using only properties of the fundamental group, and we have a homotopy equivalence between $S$ and $X$ which induces $\pi_1(S)=\pi_1(X)=\Gam$. Thus we have compatible identifications $H^4(\Shat,\ZZ)=H^4(\Xhat,\ZZ)$ and $\Alb(\Shat)=\Alb(\Xhat)$, which we arrange in the following commutative diagram:
\begin{equation*}
\renewcommand{\objectstyle}{\displaystyle}
\xymatrix@R=1pt{
  & H^4(A_{\Shat},\ZZ)\ar[r]^{a_{\Shat}^*} & H^4(\Shat,\ZZ) \\
  H^4(T,\ZZ) \ar[ur]+L^>>>>>>{i_{\Shat}^*}\ar[rd]+L^>>>>>>{i_{\Xhat}^*} &&\\
  & H^4(A_{\Xhat},\ZZ)\ar[r]^{a_{\Xhat}^*} & H^4(\Xhat,\ZZ) \ar@{=}[uu]}
\end{equation*}
Here we have defined $T=\Alb(\Xhat)\isom\Alb(\Shat)$.

Now $a_{\Xhat}^*$ is an isomorphism, and we also note that $i_{\Xhat}^*$ can not be trivial, because $A_{\Xhat}$ is a hypersurface of tridegree $(3,3,3)$ in $T$. This implies $A_{\Shat}$ is codimension $1$ in $T$, because $H^4(A_{\Shat},\ZZ)\neq0$ by commutativity. Let $[A_{\Shat}]$ and $[A_{\Xhat}]$ be the fundamental classes of
$A_{\Shat}$ and $A_{\Xhat}$ in $H^2(T,\ZZ)$ respectively. Then for any $\lambda \in H^4(T,\ZZ)$,
we have
\[(\deg a_{\Shat})[A_{\Shat}]\cdot\lambda = \alb^*(\lambda)= a_{\Xhat}^*(i_{\Xhat}^*(\lambda))=[A_{\Xhat}]\cdot\lambda.\]
By Poincar\'e duality, there exists $\lambda'\in H^4(T,\ZZ)$ such that
$[A_{\Xhat}]\cdot\lambda' = 1$, which implies that
$\deg a_{\Shat} = 1$. Thus $\Shat$ and $A_{\Shat}$ are birational, and $a_{\Shat}^*$ is an
isomorphism. Therefore, $[A_{\Shat}]\cdot \lambda = [A_{\Xhat}]\cdot \lambda$ for all
$\lambda\in H^4(T,\ZZ)$ so that $[A_{\Shat}] = [A_{\Xhat}]$, which implies $A_{\Shat}$ is a hypersurface of
tridegree $(3,3,3)$ in $T$. This proves part $(2)$ of the proposition. $\qed$
\end{pf}

To complete the proof of the theorem, we prove that the quotient of $A_{\Shat}$ by
$\Gam/\La$ is in the same moduli space as $S$. This follows from the following:
\begin{lemma} The quotient of $A_{\Shat}$ by $\Gam/\La$ is a surface with $p_g=0$, $K^2=6$ and at worst
rational double points.
\end{lemma}
\begin{pf} We compare the invariants of $A_{\Shat}$ and $\Shat$, and show
that $A_{\Shat}$ has at worst rational double points. Since
$\omega_{A_{\Shat}} = \Oh_{A_{\Shat}}(3,3,3)$, we have $K_{A_{\Shat}}^2 = 6\cdot 3^3$,
while $K_{\Shat}^2=3^3\cdot6$ because $\Shat\to S$ is an \'etale
$(\ZZ/3)^3$-cover. Note that $p_g(\Shat)= 29$ because $\chi(\Shat) = 3^3 \chi(S) = 3^3$ and $q(\Shat) = 3$.
The short exact sequence
\[0\to\Oh_T\to\Oh_T(A_{\Shat})\to\omega_{A_{\Shat}}\to0\]
gives rise to
  \begin{align*}
    0 \to H^0(T,\Oh_T) \to H^0(T,\Oh_T(A_{\Shat})) &\to H^0(A_{\Shat},\omega_{A_{\Shat}})\\
    &\to H^1(T,\Oh_T) \to H^1(T,\Oh_T(A_{\Shat})).
  \end{align*}
Since $\Oh_T(A_{\Shat})$ is very ample, we see that $|\omega_{A_{\Shat}}|$ is base point free and $H^1(T,\Oh_T(A_{\Shat})) = 0$. We also have $h^0(T,\Oh_T) = 1$, $h^1(T,\Oh_T) = 3$
and $h^0(T,\Oh_T(A_{\Shat})) = 3^3 = 27$ by K\"unneth's formula. Hence $p_g(A_{\Shat}) = 29 = p_g(\Shat)$, so $\omega_{\Shat} = \alb^*\omega_{A_{\Shat}}$ and $A_{\Shat}$ has at worst rational double points as singularities. Since the induced action of $\Gam/\La$ on $A_{\Shat}$ is free, the quotient has only rational double points, so $p_g=0$ and $K^2=6$.$\qed$
\end{pf}
Thus we have recovered the Inoue-type construction of our surface $S$ using only the fact that $S$ is homotopy equivalent to a Kulikov surface. This proves the theorem.$\qed$

\section{Degree of the bicanonical map}\label{sec!bicanonical}
Let $k_2 \colon X \to \PP^{K_X^2}$ be the bicanonical map of a Kulikov surface $X$.
Then $k_2$ is a morphism by Reider's theorem \cite{Reider}, and the image is a surface
by \cite{Xiao}. In this section we prove that the degree of $k_2$ is $1$.

Let $r$ be the degree of $k_2$, and $s$ the degree of its image in
$\PP^{K_X^2}$. Then clearly $4K_X^2 = rs$. Moreover, $r \leq 4$
because $s \geq K_X^2 - 1$, and $r\ne3,4$ by \cite{MP}, \cite{MP2} respectively, so the only possible
values for $r$ are $1$ or $2$. Now to calculate the value of $r$, we use the following proposition,
which is a natural extension of \cite{Par}, Proposition 4.1:

\begin{prop}\label{prop!bican}
Let $\psi\colon X \to Y$ be a finite abelian cover with group $G$ and branch
divisor $\Delta$, where $X$ and $Y$ are nonsingular varieties of dimension $n$. Then the
direct image of the bicanonical sheaf $\psi_*\omega_X^2$ splits into eigensheaves
 \begin{equation}
 \begin{split}\label{eqn!bican}
    (\psi_*\omega_X^2)^{\inv}&= \omega_Y^2(\Delta) \\
    (\psi_*\omega_X^2)^{(\chi)}&=\omega_Y^2
        (\sum_{g\in S_{\chi^{-1}}} \Del_g)\otimes\eL_{\chi^{-1}}
 \end{split}
 \end{equation}
where $S_\chi = \{g\in G | \frac md \chi(g)\ne m - 1\}$, $m$ is the order of $g$, $d$ is
the order of $\chi$, and $\Del_g$ is the sum of the components $\Del_i$ of $\Del$ such
that $\Fie(\Del_i)=g$.
\end{prop}
\begin{pf}

By Lemma 4.1 in \cite{Par}, locally free sheaves on $Y$ which
coincide on the complement of a codimension $\geq 2$ subset are in fact
equal on $Y$. Thus it suffices to prove the statement in a neighbourhood of
every point which lies on a single irreducible component of $\Del$.
Therefore we assume $\Del$ is irreducible and $\Del=\Del_g$ for some $g\in G$.

Let $W=X/{\Span g}$. Then $\psi$ factors through $W$ as $\psi = \rho\circ\pi$,
where $\pi\colon X\to W$ is a cyclic cover branched over $\rho^{-1}(\Del)$
and $\rho\colon W\to Y$ is unramified. Let $\eM = \Oh_W(M)$ be the
line bundle on $W$ corresponding to the cyclic cover $\pi$, so $\eM^m=\Oh_W(\rho^*\Del)$
and the formula \[\pi_*\Oh_X = \bigoplus_{\mu=0}^{m-1}\eM^{-\mu}\] decomposes $\Oh_X$ into eigensheaves.

Choose local coordinates $b, w_2,\dots, w_n$ for $W$, so that
$b$ is a local equation for $\rho^*\Del$, and let $z$ be a local generator of
$\eM^{-1}$ as an $\Oh_W$-module. Then $X$ is defined by the equation $z^m = b$,
so that $mz^{m-1}dz=db$, and we have the following local basis for $\pi_*\omega_X^2$
as an $\Oh_W$-module:
\begin{equation*}
  z^{\mu-2m+2}(db\wedge dw_2 \wedge\dots\wedge dw_n)^2,\quad\mu=0,\dots,m-1.
\end{equation*}
As in Section \ref{sec!surfaces}, we identify $\Span g^*$ with $\ZZ/m$ so that the dual character to $g$ is $1$.
Then $g$ acts on $z$ by $g\cdot z=\eps z$, where $\eps=\exp(\frac{2\pi i}{m})$,
and by considering the action of $\Span g$ on our basis, $\pi_*\omega_X^2$ splits into eigensheaves:
\begin{equation}
\begin{split}\label{eqn!pre-bican}
 \paren{\pi_*\omega_X^2}^{(0)} &= \omega_W^2(\rho^*\Del)\\
 \paren{\pi_*\omega_X^2}^{(1)} &= \omega_W^2\otimes\eM^{m-1}\\
 \paren{\pi_*\omega_X^2}^{(\mu)} &= \omega_W^2(\rho^*\Del)\otimes\eM^{m-\mu} \quad\text{for }2\leq\mu\leq m-1
\end{split}
\end{equation}

Now since $\rho$ is unramified, we have
\begin{equation*}
  \rho_*\eM^i = \bigoplus_{\frac md\chi(g)=i} \eL_\chi,
\end{equation*}
and clearly $\psi_*\omega_X^2 = \rho_*\pi_*\omega_X^2$.
Combining this with \eqref{eqn!pre-bican}, we obtain the required decomposition of $\psi_*\omega_X^2$.$\qed$
\end{pf}

In Table \ref{tab!bican} we list the eigensheaves of $\psi_*\omega^2_X$ as calculated using Proposition \ref{prop!bican}.
For example, when $\chi=(2,2)$ so that $\chi^{-1}=(1,1)$,
we have $S_{(1,1)}=\left\{\binom10,\binom01,\binom21\right\}$ and
\[\Del_{\binom10}+\Del_{\binom01}+\Del_{\binom21}=5H-2E_1-E_2-3E_3.\]
Now by equation \eqref{eqn!eigensheaf},
\[L_{(1,1)}=\frac13(\Del_{\binom10}+2\Del_{\binom11}+\Del_{\binom01})=2H-E_1-E_2,\]
so using \eqref{eqn!bican} we have
\[(\psi_*\omega_X^2)^{(2,2)}=H-E_1-E_3\]
because $K_Y=-3H+\sum E_i$.
\begin{table}[htbp]
\small
\begin{equation*}
\renewcommand{\arraystretch}{1.2}\begin{array}{|c|c|c|c|}
\hline
\chi&\sum_{S_{\chi^{-1}}}\Del_g&L_{\chi^{-1}}&(\psi_*\omega_X^2)^{(\chi)}\\
\hline
(0,0)&\Del&0&0\\
(1,0)&\Del_{\binom01}+\Del_{\binom21}&3H-E_1-2E_2-E_3&-H+E_1+E_2\\
(1,2)&\Del_{\binom01}+\Del_{\binom11}&3H-E_1-E_2-2E_3&-H+E_2+E_3\\
(1,1)&\Del_{\binom11}+\Del_{\binom21}&3H-2E_1-E_2-E_3&-H+E_1+E_3\\
(0,1)&\Del_{\binom10}&2H&-H\\
(0,2)&\Del_{\binom10}+\Del_{\binom01}+\Del_{\binom11}+\Del_{\binom21}&H&H\\
(2,2)&\Del_{\binom10}+\Del_{\binom01}+\Del_{\binom21}&2H-E_1-E_2&H-E_1-E_3\\
(2,0)&\Del_{\binom10}+\Del_{\binom01}+\Del_{\binom11}&2H-E_2-E_3&H-E_1-E_2\\
(2,1)&\Del_{\binom10}+\Del_{\binom11}+\Del_{\binom21}&2H-E_1-E_3&H-E_2-E_3\\
\hline
\end{array}
\end{equation*}
\caption{Decomposition of $\psi_*\omega_X^2$}\label{tab!bican}
\end{table}

Using the Table, we can calculate the degree of the bicanonical map of a Kulikov surface.
\begin{prop}
 The bicanonical morphism $k_2$ of a Kulikov surface $X$ is birational.
\end{prop}
\begin{pf}
It suffices to show that $k_2$ separates points in $X$.
It is clear from Table \ref{tab!bican} that the only summands of $\psi_*\omega_X^2$ with
global sections are the eigensheaves with characters $(0,2)$, $(2,2)$, $(2,0)$
and $(2,1)$ together with the invariant eigensheaf. Note that these characters
generate ${(\ZZ/3)^2}^*$.

Now choose a basis of $H^0(X,2K_X)$ comprising eigenfunctions for the
decomposition, so that $k_2$ is defined via this basis.
Then for a generic point $x$ in $X$ and any $g$ in $(\ZZ/3)^2$, if $x$ and $gx$ have the same image
under $k_2$, we have $(0,0)(g) = (0,2)(g) = (2,2)(g) = (2,0)(g) = (2,1)(g)$,
which implies $g = \binom00$. Therefore, $k_2$ separates points in each fibre of $\psi$. 

Finally we show that $k_2$ separates fibres of $\psi$. Note that 
$(\psi_*\omega_X^2)^{(0,2)}=\Oh_Y(H)$, and therefore $k_2$
separates points in $Y\smallsetminus \paren{E_1\cup E_2\cup E_3}$. It
follows that $k_2$ separates fibres of $\psi$, and so $k_2$ is birational.$\qed$
\end{pf}

\section{Bloch conjecture}\label{sec!bloch}

In this section, we use methods of \cite{IM} to show that the Bloch conjecture is verified for Kulikov surfaces. Let 
$A^0_0(S)$ denote the group of rational equivalence classes of zero cycles of degree $0$ on a surface $S$. The
conjecture is that for a surface $S$ with $p_g=0$, $A_0^0(S)$ is canonically isomorphic to $\Alb(S)$. We note that the
conjecture is proved for surfaces with Kodaira dimension $\kappa<2$ in \cite{Bloch}. Since the Kulikov surface $X$ is a surface
of general type, we must show that $A^0_0(X)=0$.

We briefly outline the approach: suppose $G$ is a finite group of automorphisms of a surface $S$. Then 
every element $g$ in $G$ induces an endomorphism $g_*\colon A^0_0(S)\to A^0_0(S)$, and this extends linearly to a
homomorphism
\[\Gam\colon\CC G\to\End(A^0_0(S)).\]
Now for a subgroup $H$ of $G$, we define 
\[z(H)=\sum_{h\in H}h.\]
Then by \cite{IM} we have $A^0_0(S/H)=0$ if and only if $\Gam(z(H))=0$. This very elegant result was
further refined in \cite{Barlow}:

\begin{lemma}\label{lem!IM} Let $H$, $H_i$, $i=1,\dots,n$ be subgroups of $G$, and let $I$ be the ideal of
$\CC G$ generated by $z(H_i)$. Suppose $A^0_0(S/H_i)=0$ for all $i$. Then if $z(H)$ is in $I$, we have $A^0_0(S/H)=0$.
\end{lemma}

Now, this lemma is used together with the Inoue construction in \cite{IM} to prove that Burniat--Inoue surfaces 
verify the Bloch conjecture. We adapt their method to Kulikov surfaces:
\begin{thm} Kulikov surfaces verify the Bloch conjecture.
\end{thm}
The proof for Kulikov surfaces is slightly more involved than in \cite{IM} because
the groups used are larger, so we present the full calculation. First we establish some notation
to streamline the algebraic manipulations which are used in the proof. Recall from Section \ref{sec!inoue} that
a Kulikov surface $X$ is a $G^2$-quotient of $\Xhat_{3,3,3}\subset E\times E\times E$, where
$G^2=\Span{g_1,g_2,g_3}\isom(\ZZ/3)^3$. Now $G^2$ is contained in the larger group $G^1\isom(\ZZ/3)^5$
of automorphisms of $\Xhat$. Following Lemma \ref{lemma!Inoue} we choose generators $\xi_1,\xi_2,\xi_3,\omega_1,\omega_2,\omega_3$
for $G^1$ subject to the relation $\xi_1\xi_2\xi_3=1$. By Lemma \ref{lemma!Inoue} and Definition \ref{def!G0G2},
 $\xi_i=\eta_i^2\eta_{i+1}$ are composite translations generating $G^0$, $\omega_i$ are rotations on the $i$th factor $E$
and $g_i=\xi_{i+1}\omega_{i+1}$.

\begin{lemma}\label{lem!rtlquotients} Let $H\subset G^1$ be one of the following collection of subgroups:
\[\Span{\omega_i,\xi_n^l\omega_j,\xi_n^m\omega_k}\text{ for }0\le l,m\le2,\]
where $\{i,j,k\}=\{1,2,3\}$ and $n=i$ or $i-1$. Additionally we may choose
$H=\Span{\xi_3\omega_1,\xi_2\omega_2,\xi_2\omega_3}$. Then the surface $\Xhat/H$ is rational.
\end{lemma}
\begin{rmk}\rm
This is not an exhaustive list of subgroups $H$ which give rise to a rational quotient, but it is sufficient
to prove the Theorem. The additional subgroup $\Span{\xi_3\omega_1,\xi_2\omega_2,\xi_2\omega_3}$ is necessary
for the proof, and this is precisely what makes the calculation harder than the elementary one in \cite{IM}.
\end{rmk}
\begin{pflemma}
Recall the definition of $\Xhat$ as a $G^1$-cover of the plane branched in the Kulikov line configuration
from Section \ref{sec!kulikovsurf}. In particular, by Definition \ref{def!G1}, the columns $\del_1$, $\del_2$,
$\del_3$, $\omega_1$, $\omega_2$, $\omega_3$ of matrix \eqref{eqn!kulikovFie} generate $G^1$. Then by
Definition \ref{def!G0G2} we have
\[\xi_i=\del_{i-1}^2\omega_i\omega_{i+1},\]
and since $\xi_1\xi_2\xi_3=1$, we also have
\[\del_1\del_2\del_3=(\omega_1\omega_2\omega_3)^2.\]

Now suppose $H$ is a subgroup of $G^1$. The quotient $\Xhat/H$ is a ${G^1}/H$-cover of $\PP^2$ branched in the Kulikov line
configuration of Figure \ref{fig!configs}. Thus we may use the techniques of abelian covers to prove that $\Xhat/H$ is rational.
This involves a series of repetitive calculations, and various cases are related to one another by symmetry. As an illustration,
we calculate a typical quotient:

Let $H=\Span{\omega_1,\xi_3^2\omega_2,\xi_3\omega_3}$, then ${G^1}/H\isom(\ZZ/3)^2$ is generated by the classes $\coset{\del_1}$ and $\coset{\omega_3}$.
We list the class of each generator of $G^1$ under the quotient:
\[\coset{\del_2}\equiv\coset{\omega_3^2},\quad\coset{\del_3}\equiv\coset{\del_1^2\omega_3},
\quad\coset{\omega_1}\equiv1,\quad\coset{\omega_2}\equiv\coset{\omega_3^2}.\]
Viewing $\Xhat/H$ as a $(\ZZ/3)^2$-cover of $\PP^2$ branched in the Kulikov line configuration,
we must blow up $P_2$ and $P_3$ to remove singular points on the cover. Moreover, the point of intersection of $\Del_5$
and $\Del_6$ fails the nonsingularity condition of equation \eqref{eqn!nonsingular}. Thus we blow up this point, introducing
another $-1$-curve which we call $E$. Having done this, we note that by equation \eqref{eqn!exceptional}, the induced cover
is unramified over $E_2$ and $E$. Finally, we contract the strict transforms of $\Del_1$ and $\Del_5$. We are left with a
nonsingular $(\ZZ/3)^2$-cover of $\PP^1\times\PP^1$ branched in four distinct lines, two in each ruling. This is a rational surface.$\qed$
\end{pflemma}
\begin{pfthm}
Let $I$ be the ideal generated by $z(H)$ for all $H$ listed in Lemma \ref{lem!rtlquotients}. Then the quotients $\Xhat/H$ are
rational surfaces, so by \cite{Bloch}, $A^0_0(\Xhat/H)=0$. Thus by Lemma \ref{lem!IM}, it is
sufficient to prove that $z(G^2)$ is an element of $I$. Now, the verification of the Bloch conjecture for Burniat--Inoue
surfaces \cite{IM} is a series of elementary polynomial manipulations. Unfortunately, the corresponding manipulations do
not suffice for the Kulikov surface.

Instead, we prove that $z(G^2)$ is in $I$ using the Magma computer algebra script \cite{Magma} below:

\begin{verbatim}
Q:=Rationals();
RR<w1,w2,w3,xi1,xi2,xi3>:=PolynomialRing(Q,6);

function z(a)
   return (1+a[1]+a[1]^2)*(1+a[2]+a[2]^2)*(1+a[3]+a[3]^2);
end function;

ListH:=[[w1,xi1^i*w2,xi1^j*w3]:i in [0..2],j in [0..2]] cat
       [[w1,xi3^i*w2,xi3^j*w3]:i in [0..2],j in [0..2]] cat
       [[w2,xi1^i*w3,xi1^j*w1]:i in [0..2],j in [0..2]] cat
       [[w2,xi2^i*w3,xi2^j*w1]:i in [0..2],j in [0..2]] cat
       [[w3,xi2^i*w1,xi2^j*w2]:i in [0..2],j in [0..2]] cat
       [[w3,xi3^i*w1,xi3^j*w2]:i in [0..2],j in [0..2]];

Append(~ListH,[xi3*w1,xi2*w2,xi2*w3]); // extra generator

I:=ideal<RR|[xi1^3-1,xi2^3-1,xi3^3-1,
             w1^3-1,w2^3-1,w3^3-1,xi1*xi2*xi3-1]
            cat [z(H):H in ListH]>;

zG2:=z([w1*xi1,w2*xi2,w3*xi3]);

zG2 in I; // Result: true
\end{verbatim}
$\qed$
\end{pfthm}

\vspace{0.2cm}
\textbf{Authors' address}\\
\small{
Lehrstuhl Mathematik VIII,\\Mathematisches Institut der Universit\"at Bayreuth,\\NW II, Universit\"atstr. 30,\\95447 Bayreuth}
\begin{verbatim}
mariocto@gmail.com
stephen.coughlan@uni-bayreuth.de
\end{verbatim}

\end{document}